\documentclass[11pt]{amsart}
\usepackage{amssymb}

\usepackage{palatino}
\input amssym.def

\usepackage{amsmath, amsfonts}
\usepackage{amssymb}
\usepackage{amscd}
\usepackage[mathscr]{eucal}
\usepackage{palatino}
\newfont{\cyrr}{wncyr10}
\newcommand{\bO}{{\mathcal O}}
\renewcommand{\H}{{\mathcal H}}

\renewcommand{\a}{{\mathfrak a}}
\renewcommand{\b}{{\mathfrak b}}
\newcommand{\gC}{{\mathfrak C}}
\newcommand{\N}{{\mathbb N}}
\newcommand{\Z}{{\mathbb Z}}
\newcommand{\Q}{{\mathbb Q}}
\newcommand{\R}{{\mathbb R}}
\newcommand{\C}{{\mathbb C}}

\newtheorem{thm}{Theorem}

\newtheorem{lem}[thm]{Lemma}

\newtheorem{rmk}{Remark}[section]

\begin{document}

\title[]{ Around a question of Baker }

\author{Neelam Kandhil and Purusottam Rath}

\address{Neelam Kandhil \\ \newline
The Institute of Mathematical Sciences, A CI of Homi Bhabha National Institute, 
CIT Campus, Taramani, Chennai 600 113, India.}
\email{neelam@imsc.res.in} 

\address{Purusottam Rath\\ \newline
Chennai Mathematical Institute,  Sipcot IT Park, Siruseri, Kelambakkam 603103,
India.}

\email{rath@cmi.ac.in} 

\subjclass[2020]{11J72, 11J86, 11R18}

\keywords{Dirichlet $L$ functions, Linear Independence over Number fields}

\maketitle

\begin{abstract} 
For any positive integer $q$,
it is a question of Baker whether the  numbers  $L(1, \chi)$,
where $\chi$ runs over the non-trivial characters mod $q$,
are linearly independent over $\Q$. The question is 
answered in affirmative for primes, but
is unknown for an arbitrary modulus.
In this expository note,
we give an overview of  the origin and history  of this question as well as the state-of-the-art.
We also give an account of the mathematical ideas that enter into
this theme as well as elucidate the obstructions that preclude us from 
answering the question for arbitrary modulus. We also describe 
a number of  generalisations and extensions of this question.

\end{abstract}

\section{Introduction}

The central characters of this article are the Dirichlet characters.
These  are  one dimensional Galois representations 
of Cyclotomic extensions. More concretely,
for an integer $n>1$, a Dirichlet character $\chi$ is  simply a homomorphism from
the group  $(\Z/n\Z)^{\times}$
of co-prime residue classes mod $n$ 
 to the multiplicative group $\C^{\times}$. By assigning the value zero at the  other classes mod $n$, we can extend
$\chi$  to a function from $\Z$ to $\C$ which is completely
multiplicative and periodic with  period $n$.

Each such character $\chi$ gives rise to a Dirichlet series
$$
L(s, \chi) = \sum_{n=1}^{\infty} \frac{\chi(n)}{n^s},
$$
where the series is absolutely convergent in the region $\Re(s) > 1$.
Furthermore, since $\chi$ is completely multiplicative, one has the 
Euler product representation
$$
\sum_{n=1}^{\infty} \frac{\chi(n)}{n^s} = \prod_{p ~~ {\rm prime}}
\left( 1 - \chi(p) p^{-s}\right)^{-1}
$$
which is a consequence of prime factorisation of integers.
When $\chi$ is a non-trivial character,  we know that $L(s,\chi)$
extends to an entire function and  we have
$$
L(1, \chi) = \sum_{n=1}^{\infty} \frac{\chi(n)}{n}
$$

and it is these complex numbers which are the centre of our focus. From
now on, the Dirichlet characters are assumed to be non-trivial, unless
stated otherwise.

 A celebrated
result of Dirichlet asserts that  $L(1, \chi)$ is non-zero. This work of Dirichlet  laid
 the foundation of Analytic Number theory. This also ushered
in the application of Character theory into the realms of Number theory,
a theme which has now spectacularly morphed into the enigmatic
interplay between Harmonic Analysis and Arithmetic of Galois representations.
Furthermore, these special values have deep arithmetic content.
For instance, for any quadratic number field $K$, one has
a quadratic Dirichlet character $\chi$ associated to $K$ and
$L(1,\chi)$ subsumes deep arithmetic data like
the class number and Regulator of $K$.

Let us now come to the focal point of this note. Since the $L(1,\chi)$'s are non-zero,
it is natural to ask about the algebraic nature of these complex numbers.
While the non-vanishing was established in 1837, it took about
130 years more to settle the nature of these numbers. The seminal work
of Baker in mid 1960's established that these numbers are transcendental.

For a positive integer $q$,  each of the $\varphi(q) - 1$ non-trivial Dirichlet characters
give rise  to seemingly  unrelated   transcendental numbers.  But these
numbers do not  quite live in complete isolation from each other.
Let $K = \Q(\zeta_q)$ be the $q$-th Cyclotomic field and $\zeta_K(s)$
be its Dedekind zeta function.
The product of these $L(1,\chi)$ values gives the residue of  
$\zeta_K(s)$ at $s=1$. So it is natural to wonder about any
relation, linear or algebraic, existing between these mysterious numbers.

In 1973, Baker, Birch and Wirsing \cite{bbw} in a beautiful work
proved that for a prime $p$, the  numbers  $L(1, \chi)$
where $\chi$ runs over the non-trivial characters mod $p$,
are linearly independent over $\Q$.
Baker, in his book Transcendental number Theory (\cite{baker}, p.48), stated that it would be
of interest to know if this is true for an arbitrary modulus $q$.
This remains unanswered till now and is the raison d'etre of our note.

We attempt to give a broad account of the history and the  state-of-the-art
of this question of  Baker. We highlight the mathematical tools
and techniques that enter into this circle of questions.
We also describe  a number of  generalisations and extensions of this question, namely
extensions to Number fields, to class group L-functions as well
as specialisations at larger  integers.
There are some essential ingredients common to each of the above themes,
viz, Galois
action on linear forms of certain periods, Transcendence theory of linear forms in logarithms,
non-vanishing of L-values from  Arithmetic
and finally Dedekind-Frobenius determinants from Linear algebra. We shall try
to illustrate  the commonality in the   above circle of questions
as well as the issues intrinsic to each of these separate themes.

Let us end this section  by noting that the complex numbers $L(1,\chi)$'s 
 are examples of {\it Periods}. A  period  as defined by Kontsevich and Zagier \cite{KZ}
 is a complex number whose real and imaginary parts
are values of absolutely convergent integrals of rational functions
with rational coefficients over domains in ${\R}^n$
given by polynomial inequalities with rational coefficients.
Clearly, all algebraic numbers are periods.  An important class of
periods is supplied by the special values of the Riemann
zeta function.
For example, we have for $k \ge 2$
$$\zeta(k) = \int_{1 > t_1 > \cdots > t_k >0}
{dt_1 \over t_1} \cdots {dt_{k-1} \over t_{k-1}} {dt_k \over 1 - t_k}, $$
as is easily verified by direct integration.  Also $\pi$ is a period, it is expected
that $e$ is not a period. 

The set of periods  is countable and a $\bar{\Q}$-algebra of infinite dimension.
Since logarithms of algebraic numbers are periods, we shall  see 
that each $L(1,\chi)$ is indeed a period.
We recommend the original delightful  article of Kontsevich and Zagier \cite{KZ} as well as the
account of Waldschmidt \cite{waldperiods} for  further details.
We also heartily
recommend the comprehensive book by Huber and M$\ddot{u}$ller-Stach \cite{huber}.

\section{A question of Chowla and the remarkable result \\
of Baker, Birch and Wirsing}

In a lecture at the Stony Brook conference on number theory in
 1969, Sarvadaman Chowla posed the  question whether there exists a non-zero rational-valued arithmetic function $f$, periodic with
prime period $p$ such that $\sum_{n=1}^\infty {f(n) \over n} = 0.$

 In 1973, Baker, Birch and Wirsing \cite{bbw} 
answered this question in the following theorem:
\index{Baker, A.} \index{Birch, B.} \index{Wirsing, E.}
\begin{thm}  If $f$ is a non-zero function defined on the
integers with algebraic values and period $q$ such that
$f(n)=0$ whenever $1< (n,q)<q$ and the $q$-th cyclotomic polynomial
is irreducible over the field $F = \Q(f(1), ..., f(q))$, then
$$\sum_{n=1}^\infty {f(n) \over n} \neq 0. $$
\end{thm}

In particular, if $f$ is rational valued, the second condition
holds trivially.  If $q$ is prime, then the first condition is
vacuous.  Thus, the theorem resolves Chowla's question.

\smallskip

The above theorem of  Baker, Birch and Wirsing is remarkable as
it  brought out a new aspect of Transcendence theory hitherto undiscovered.
More often than not,  non-vanishing is a major obstacle in Transcendence and
typically non-vanishing follows from Arithmetic considerations. For instance,
 transcendence of $L(1,\chi)$ follows only when its non-vanishing is ensured.
But the ideas in the proof of Baker-Birch-Wirsing indicated that 
the tables can be turned and transcendence can ensure non-vanishing.
This perspective  has been
exploited quite fruitfully in recent times. For instance, Kumar Murty and Ram Murty \cite{ram-kumar-jnt}
have used Transcendence theory to prove the non-vanishing
of $L(1,\chi) (\chi \neq 1)$ for even characters, an illustrative  instance
of Transcendence theory returning the favours to Arithmetic. 
Recall that a character $\chi$ is even or odd according as $\chi(-1)$ is $1$ or $-1$.

\smallskip

Let us now highlight the main ingredients  in the proof of  Baker-Birch-Wirsing Theorem.

\smallskip
\noindent
1.  The series  $\sum_{n=1}^{\infty} \frac{f(n)}{n}$ converges
if and only if  $\sum_{a =1} ^ q  f(a) = 0$. Once this is ensured,
one derives that 
$$
{\rm L}(1,f) = \sum_{n=1}^{\infty} \frac{f(n)}{n}
=  \frac{-1}{q}  \sum_{(a,q)=1} f(a)  (\psi(a/q) + \gamma).
$$
Here  $\psi$ is the digamma function which  shows up since it is  the constant term  of the Hurwitz zeta function
around  $s=1.$ We will need to come back 
to the Hurwitz zeta function in the last section. 
But the digamma function  regrettably is a rather difficult
function to handle, for instance, $-\psi(1)$ is the enigmatic
Euler's constant  $\gamma$.

\smallskip
\noindent
2. Since $f$ is periodic, we can Fourier analyse and go to  the dual set up.
Let 
$$
\hat{f}(n) = { 1 \over q} \sum_{m=1}^{q} f(m) e^{-2\pi imn/q}
$$
be the  Fourier transform of $f$. By inverse Fourier transform and 
some functional manipulation, one is led to the following new identity
$$
{\rm L}(1,f) =  -\sum_{a=1}^{q-1} \widehat{f}(a) \log(1- \zeta_q^a).
$$
While we no longer have the digamma function, the caveat
 now is the added complexity of the 
  sequence of Fourier coefficients $\{\hat{f}(n)\}$ as they need not
lie in the field generated by the original  sequence $\{f(n)\}$. But the redeeming
feature is that the  digamma function is replaced by 
logarithms of algebraic numbers, setting the tone
for the entry of Baker's seminal theory.
Without further  ado, let us state  Baker's theorem (see \cite{baker}, for instance) which is a  pivotal
ingredient in our context.

\begin{thm} \label{baker} {\rm (Baker)}  If $\alpha_1, ..., \alpha_m$
are non-zero algebraic numbers such that 
$\log \alpha_1,  ..., \log \alpha_m$ are  linearly independent over $\Q$,
then $1, \log \alpha_1, ..., \log \alpha_m$ are linearly independent
over $\overline{\Q}$.  
\end{thm}

\smallskip
\noindent
3. Baker's theorem allows one  to show that  vanishing of $L(1,f)$
ensures the vanishing of the "conjugate"  L-values $L(1, f^{\sigma})$
for $\sigma$ in the Galois group  $G$ of the extension ${F}(\zeta_q)/ {F}$.
Let us be more concrete. We note
that $G$ is isomorphic to the group $(\Z/q\Z)^{\times}$. For $(h,q) = 1$, let $\sigma_h \in G$
be such that
$
\sigma_h(\zeta_q) = \zeta_q^h.
$
Define $ f^{\sigma_h} = f_h(n) := f(nh^{-1})$ for $(h,q)=1$.
Then, Baker, Birch and Wirsing  beautifully exploited  Baker's theorem
to show that the vanishing of $L(1,f)$ results in
$$
{\rm L}(1,f_h) = \displaystyle \sum_{n=1}^{\infty} \frac{f_h(n)}{n} 
= 0
$$
for all $(h,q) = 1 $.

\smallskip
\noindent
4. We then reverse the  Fourier process and  come back to
the original set up
$$
{\rm L}(1,f_h) = \frac{-1}{q} \sum_{(a,q) = 1 } f_h(a)
(\psi(a/q) + \gamma) = 0
$$
leading to the following  system of identities (for each $h$  co-prime to $q$)
$$ \sum_{(a,q)=1} f(a) (\psi({ah}/q) + \gamma) = 0.
$$

\smallskip
\noindent
5. We now notice that the matrix $
{\rm A} : =  \left( \psi({ah}/q) + \gamma \right)_{(ah,q) = 1}
$
associated to the system of identities obtained in the previous step
is a Dedekind-Frobenius matrix on the group   H  =  $(\Z/q\Z)^{\times}$ and its determinant 
(up to a sign) is given by  
$$
\prod_{\chi \in \widehat{\rm H}} \left(\sum_{h \in {\rm H}} \chi(h) (\psi(h/q) + \gamma) \right).
$$
If we show that the matrix A is invertible, then $f$ vanishes everywhere and we are done.

\smallskip
\noindent
6. This is where  "Arithmetic"  enters the picture. The non-vanishing
of determinant of $A$ is ensured by the non-vanishing of  ${\rm L}(1, \chi)$ for non-trivial $\chi$
while the monotonicity of $\psi$ takes care of the  trivial character.
 This completes the proof of the theorem.

\smallskip
\noindent
We mention in passing that in \cite{gmr}, a generalization of the above theorem
has been derived.

{\rmk
{\it
In 1966,  Lang \cite{Lang1} proved a multi-dimensional generalisation
of the classical theorem of Schneider, motivated by  a question of
Cartier about the analogue  of transcendence of $e^{\alpha}, \alpha \in \bar{\Q}^{\times}$
for arbitrary group varieties. Lang did answer this question in affirmative. 
In retrospect,  as evinced by a later work of Bertrand and Masser,
we see that Baker's Theorem  in the form stated above
could have been proved by  Lang  in 1966 
 building on the Galois conjugate idea in Baker-Birch-Wirsing
Theorem and his generalisation of Schneider's theorem.
 But perhaps it is fortuitous that Lang did not prove this theorem.
It is because the approach of Baker is different who not only proved  the above
theorem, but also obtained lower bounds for linear forms in logarithms of algebraic numbers.
These lower bounds are of seminal importance in the study of  diophantine equations and
form a subject of its own. For instance, these lower bounds  allowed Baker
to  classify all imaginary quadratic fields
with class number one, a venerable theme in Number theory set in motion
by Gauss. 

As indicated above, Bertrand and Masser \cite{BM} gave a new proof
of Baker's theorem by Galois action on  linear forms. They also exploited these ideas 
to prove an elliptic analog of Baker's
theorem.  Let  us state this result.
For a Weierstrass $\wp$-function with algebraic invariants $g_2$ and $g_3$ and field 
of endomorphisms $k$, the set
$
\mathcal{L} = \{ \alpha \in \C  ~: ~  \wp(\alpha) \in \overline{\Q} \cup \{\infty\} \}
$ is the two-dimensional analogue of logarithms of algebraic numbers.
Bertrand and Masser proved that if 
 $u_1, \cdots, u_n \in  \mathcal{L}$ are linearly independent over $\Q$, then
$1, u_1, \cdots, u_n$
are linearly independent
over $\overline{\Q}$.  }}

\bigskip
\noindent

\section{Settling for prime modulus and over $\overline{\Q}$}

Let us begin by noting that the above theorem of Baker, Birch and Wirsing
settles the question of Baker for prime modulus. It is because the 
values taken by the Dirichlet characters mod $p$ lie in the field
$\Q(\zeta_{(p-1)})$ which is linearly disjoint with the $p$-th
cyclotomic field. So all the hypotheses of Baker-Birch-Wirsing
theorem are satisfied and hence for an odd prime $p$,
the numbers  $L(1, \chi)$
where $\chi$ runs over the non-trivial characters mod $p$,
are indeed linearly independent over $\Q$. In a recent work \cite{GK}, the 
above result has been extended to any arbitrary family of moduli.

Let us now consider the 
possible $\bar{\Q}$-linear relations between these
values of $L(1,\chi)$ as $\chi$ ranges
over all non-trivial Dirichlet characters mod $q$ with $q>2$.
It is remarkable that over $\bar{\Q}$, we have a complete
answer for all moduli  $q$. This follows  from a natural
extension of  the works of Ram and Kumar Murty \cite{ram-kumar-jnt}.
Let us give a summary of their work.
One of the crucial results   in their work is the following:
\begin{thm} For any integer $q > 2$, the numbers  $L(1,\chi)$ as $\chi $ ranges over 
 non-trivial even characters mod $q$ are linearly independent over $\bar{\Q}$.
\end{thm}

{\rmk
This in particular furnishes a new proof of non vanishing of $L(1,\chi)$
for even non-trivial characters  by transcendental means.
The possibility of  such an approach  could not have been envisaged, but for
the work of Baker-Birch-Wirsing. Furthermore  this result
shows that the dimension of  space generated by $L(1, \chi)$
for even characters remains the same over any number field. As we shall see
a little later, this  is a luxury which is not at all afforded to us for $L(k, \chi)$ with $k > 1$.
}

\smallskip
\noindent

Let us briefly indicate the main points in the proof of this result.
The new ingredient in the proof of the above theorem
is the properties of a set of 
real multiplicatively 
independent units in the cyclotomic field discovered by
K. Ramachandra \index{Ramachandra, K.} (see Theorem 8.3 on page 147 of
\cite{wash} as well as \cite{ram2}). These marvellous units allowed 
Ram  and Kumar Murty to work with new expressions of  $L(1,\chi)$
for even characters in terms the logarithms of positive real numbers.
Thereafter, they appeal to Baker's theorem and its variants 
leading to a system of equations involving characters of finite groups.
Finally, they prove a variant of Artin's theorem
on linear independence of irreducible characters 
which establishes the desired linear independence over $\bar{\Q}$. 
 Let us state this elegant group theoretic  result for the sake of completion.
Let $G$ be a finite group.
Suppose that $\sum_{\chi\neq 1} \chi(g)u_\chi = 0 $
for all $g \neq 1$ and all irreducible characters $\chi \neq 1$ of $G$.
Then $u_\chi = 0$ for all $\chi \neq 1$.

So the space generated by the even characters is now well and truly done.
Let us now consider the odd $L(1, \chi)$ values. We note that for any odd 
Dirichlet character $\chi$, $L(1,\chi)$ is  a non-zero algebraic multiple of $\pi$.
This follows from the expressions we indicated in the previous section
for $L(1,f)$  after applying Fourier transform of $f$. Therefore, the space generated 
by the odd $L(1, \chi)$ values is one dimensional over $\bar{\Q}$.

But do these subspaces intersect? Here we come to the following  pretty
application of Baker's theory  which has been proved
in \cite{msaradha} and has a large number
of applications in various  different set ups.

\begin{lem}  Let $\alpha_1, \alpha_2, ..., \alpha_n$ be 
positive algebraic numbers.
If $c_0, c_1, ..., c_n$ are 
algebraic numbers with $c_0\neq 0$, then
$$c_0 \pi + \sum_{j=1}^n c_j \log \alpha_j $$
is a transcendental number and hence non-zero.
\end{lem}

The above  leads to the following result since  for an even character
$\chi \ne 1$, the number $L(1,\chi)$ is a linear form in logarithms of positive real algebraic numbers.

\begin{thm}
The $\overline{\Q}$-space generated
by the (non-trivial) even $L(1,\chi)$ values is linearly disjoint
from the space generated by the
odd $L(1,\chi)$ values.
\end{thm}

Consequently, one has the following satisfying result.

\begin{thm}  For any integer $ q>2$, the $\bar{\Q}$-vector space generated
by the values $L(1,\chi)$ as $\chi$ ranges
over the non-trivial Dirichlet characters (mod $q$)
has dimension $\varphi(q)/2$.
\end{thm}

Let us end this section with specifying what exactly is the issue 
which hinders us from settling Baker's question for an arbitrary modulus.
As we noted in the proof of Baker-Birch-Wirsing therem,  the
central part of the proof was to show that the vanishing of
$L(1,f)$ ensues the vanishing of the Galois conjugates $L(1,f^{\sigma})$.
A careful look in the  proof of this part reveals that it was important
that the Galois elements $\sigma$ kept the coefficients 
$\{{f}(n)\}$ unchanged. Let us give some more indication
why and when this comes up.
In the course of the proof, 
Baker's theory eventually leads to a family of expressions (indexed by $b$) of the form
$$\sum_{a=1}^{q-1} \widehat{f}(a) {\rm r}_{ab} = 0,$$
where the numbers ${\rm r}_{ab}$ lie in the field of definition $F$ of $f$. 
Then for any automorphism $\sigma \in {\rm Gal}({F}(\zeta_q)/ {F})$,
we need to act $\sigma$ on each of this identities. 
The  condition of irreducibility of the $q$-th cyclotomic polynomial over $F$ ensures that
$$
\sum_{a=1}^{q-1} \sigma(\widehat{f}(a)) {\rm r}_{ab} = 0, $$
which  then leads to 
$$
 \displaystyle \sum_{a=1}^{q-1} \sigma(\widehat{f}(a)) \log(1 - \zeta_q^{a}) =0
$$
which is what we desire.
In general, if the  automorphisms do act nontrivially on the sequence $\{f(n)\}$,
the issue becomes involved, some instances
of which we shall see in the next section.

On the other hand, the space generated by the even $L(1,\chi)$
values present no problem, thanks to  the result of Ram and Kumar Murty described 
in this section. So it is only the linear independence of
odd $L(1, \chi)$ values which remain unresolved for an arbitrary modulus.
But since these are algebraic multiples of $\pi$,  perhaps transcendence
theory   has no more role to play.
This parity conundrum will show up again a bit later when we work with $L(k, \chi)$
with $k>1$.

\section{Extension to Number fields}

In this section, we consider the extension of Baker's question to arbitrary 
number fields. Since the question is open for arbitrary modulus over
$\Q$, it is prudent to consider  the number field extension
only for prime modulus for now. This constitutes
the ethos of a recent work \cite{rathbhar}.

Let $K$ be a number field and $p$ be an odd prime.
Let us consider  the $K$-vector space in $\C$ generated by the
 $L(1, \chi)$ values for non-trivial characters $\chi$ modulo $p$. Let $d(K,p)$ denote
its dimension.
In view of the discussions in the previous section, we have the following bounds:
$$
\frac{p-1}{2} \leq d(K,p) \leq {p-2}.
$$

When $K =\Q$, the upper bound is attained. On the other hand, when
$K = \Q(\zeta_p, \zeta_{p-1})$, the lower bound is attained.
Therefore, one can ask the  following question: Which 
  numbers in the interval $$\left(\frac{p-1}{2}~~,~~ p-2\right)$$
 can be  equal to  $d(K,p)$ as $K$ runs over all number fields?
This is not known. From now onwards all primes
are at least $7$.

The next question is to ask whether for any prime $p>5$, there is 
a number field such that  
$$
\frac{p-1}{2} < d(K,p) < {p-2}.
$$

This question is answered in the affirmative in \cite{rathbhar}.
The initial strategy in \cite{rathbhar} is to look for primes
of specific type which may be more amenable 
to work with. The  family of primes 
which seem more tractable in this context are the 
Sophie Germain Primes.
A  prime $p$ is called a {\it Sophie Germain prime} if $2p+1$ is also a prime.
It  is a folklore  conjecture that there are infinitely many Sophie Germain primes.
Following theorem is proved in \cite{rathbhar}.
\begin{thm}\label{1}
Let $p > 5$ be an odd prime. Then 
there exists  a number field
$K$ such that 
$$
\frac{p-1}{2} < d(K,p) < {p-2}.
$$
\end{thm}

This theorem is proved in two steps.
In the first step,  one considers  primes that are not 
Sophie Germain where one investigates  arithmetic of  number
fields $K$ with $\Q(\zeta_{p-1}) \subset K \subset \Q(\zeta_{p-1}, \zeta_{p})$.
In the second step, one  proves the result for Sophie Germain primes by working with fields 
$K$ such that $\Q(\zeta_{p}) \subset K \subset \Q(\zeta_{p-1}, \zeta_{p})$.

Consequently for any $p> 7$, at least one number in the interval
$$\left(\frac{p-1}{2}~~,~~ p-2\right)$$
is realised as the dimension of space of $L(1,\chi)$ values
over some number field. Let  $b(p)$ count the
numbers in the above  interval that can be realised as this.
More precisely, for a prime $p>5$, let
 $$
b(p) := \left\vert\left\{n~{\Big{|}} \phantom{a}  \frac{p-1}{2} <  n < {p-2} ~~~{\rm and} ~~~~~d(K,p) =n
~~~~~~{\rm for ~~some ~~number~~ field} ~~~~ K \right\}\right\vert.$$

The above theorem implies that  $b(p) > 0$  for every prime $p>5 $ and hence
we have,
$$
1 \leq b(p) \leq \frac{p-3}{2}.
$$
Then in \cite{rathbhar}, the following is proved.
\begin{thm} \label{3} The sequence $\{b(p)\}$ satisfies
$$
\limsup_{p \to \infty} ~~b(p) = \infty.$$
 \end{thm}

Thus the sequence $\{b(p)\}$ is unbounded. One can ask
about its growth. In \cite{rathbhar}, the following
Omega result is established.

\begin{thm} \label{4} There exists a constant $c > 0$ such that
$$
b(p) >  {\rm exp} \left( \frac{ c \log p}{\log\log p}\right)
$$
for infinitely many primes $p$.
 \end{thm}

For the proof of these results, one has to work with families of number fields
 for which  the methods and approaches of the earlier sections no longer 
 work and therefore we shall not dwell further. We shall just indicate
 one of the ingredients, a  folklore  result of  Linnik (\cite{linnik1}, \cite{linnik2})
 which constitutes a celebrated theme  in Analytic number theory
 with far-reaching implications.  
\begin{thm} 
Let $a, n$ be two positive integers with $(a,n) = 1, n \geq 2$.
Let $p(a,n)$ denote the least prime $p$ such that $p \equiv a~ ({\rm mod}~ n)$.
There exists absolute positive constants $C$ and $L$ such that
$$
p(a,n) <  C n^L.
$$
\end{thm}

The constant $L$ is known as  {\it ''Linnik's constant''}. 
It is conjectured that $p(a,n) <   n^2.$ The best known value for $L$ is due to 
Xylouris \cite{xy} who proves that $L$ can be taken to be
$5.18$. On the other hand,
Lamzouri, Li and Soundararajan \cite{sound} have shown, under Generalised Riemann Hypothesis, that
$p(a,n) \leq \varphi(n)^2 \log^{2} n$ for all $n>3$. 

\section{Analogous question for class group L-functions}

In this section, we consider the analog of the 
question of Baker for class group $L$-functions.
This has been carried out by Ram and Kumar Murty \cite{murty-murty}.

Let $K$ be a number field. Let $\bO_K$ be its ring of integers and $\H_K$
be its  ideal class group. It is this finite group on which
our functions will act.

Let $f$ be
a complex-valued function of the ideal class 
group $\H_K$ of $K$. For such an $f$, we
 consider the Dirichlet series for $\Re(s) > 1$
\begin{equation*}\label{series-2}
L(s,f) := \sum_\a  {f(\a  )\over \N(\a )^s},
\end{equation*}
where the summation is over non-zero ideals $\a$ of  $\bO_K$.  
If  $f\equiv1$, $L(s,f)$ is simply the Dedekind zeta function
of $K$. We hope to study the numbers $L(1, f)$
as and when they exist.  A necessary and sufficient condition of existence is  given by the following (see \cite{rath-murty}
as well as \cite{lang1}):

\begin{thm}\label{continuation}
$L(s,f)$ extends analytically for all $s \in \C$ except
possibly at $s=1$ where it may have a simple pole with residue a non-zero multiple of
$$\rho_f : = \sum_{\a \in \H_K} f(\a). $$
Consequently, the series $
 \sum_\a  {f(\a  )\over \N(\a )}$ converges and equal to $L(1,f)$   if
and only if $\rho_f = 0$.
\end{thm}

We want  to investigate the values  $L(1,f)$ 
when $K$ is  imaginary quadratic and when $f$ takes algebraic values,
 in particular  the values $L(1,\chi)$ when $\chi$
runs over  ideal class characters. We note that  complex conjugation acts on the group of ideal class characters
and  $L(1,\chi)=L(1,\overline{\chi})$
for any ideal class character $\chi$. Let  $\H_K^*$ denote 
a set of orbit representatives under this action. Here is a pretty result proved in  \cite{murty-murty}.

\begin{thm} \label{indep} Let $K$ be an imaginary quadratic field and
$\H_K$ its ideal class group.  The values $L(1,\chi)$ as
$\chi$ ranges over the non-trivial characters of $\H_K^*$ and $\pi$ 
are linearly independent over $\overline{\Q}$.
\end{thm}

We add that unlike the case of  Dirichler characters, the above does not
prove the transcendence of $L(1,\chi)$. However it does
prove that at most one  of the values $L(1,\chi)$, as $\chi$
ranges over the non-trivial characters of $\H_K^*$,
can be algebraic.

Recall that Baker's question for primes is a consequence of the Baker-Birch-Wirsing Theorem.
Here is the analogue of Baker-Birch-Wirsing of which the above theorem is an immediate consequence.

\begin{thm}  \label{bbw-analogue}  Let $K$ be an imaginary
quadratic field and $f:\H_K \to \bar{\Q}$ be not identically zero.
Suppose that $\rho_f =0$.  
Then, $L(1,f) \neq 0$ unless $f(\gC)+f(\gC^{-1})=0$ for
every ideal class $\gC \in \H_K$.\end{thm}

So we need to prove the above theorem.
Let us indicate the salient features in the proof  the above,
indicating the commonality
with that of Baker-Birch-Wirsing as well as the new
ingredients intrinsic to this set up.
One needs to  use Kronecker's limit
formula as discussed in the works of Siegel \cite{siegel},
Ramachandra \cite{ramachandra} and Lang \cite{langel}.
We shall need Baker's theorem from Transcendence theory  as well as  Chebotarev density theorem
from Algebraic number theory.
Finally, we shall need  a few results from  Theory
of Complex Multiplication. However, we  can
only take a cursory glance  into this delightful realm 
and  shall enthusiastically direct the interested reader to
the original work \cite{murty-murty}  (and to  \cite{cox} and \cite{langel}).

We begin with  the celebrated discriminant functions $\Delta(z)$:
$$\Delta(z) = (2\pi)^{12} ~q\prod_{n=1}^\infty (1-q^n)^{24} = (2\pi)^{12}\eta(z)^{24}, \qquad q=e^{2\pi iz} $$
where $\eta^{24}$ is the ubiquitous Ramanujan cusp form.
As before, let $K$ be an imaginary quadratic field and 
 $\b$ be an ideal of $\bO_K$.  If $[\beta_1, \beta_2]$ is
an integral basis of $\b$ with $\Im(\beta_2/\beta_1)>0$, we define
$$g(\b) = (2\pi)^{-12}(\N(\b))^6 |\Delta(\beta_1,\beta_2)|,$$
where $$\Delta(\omega_1,\omega_2) := \omega_1^{-12}\Delta\left(\frac{\omega_2}{\omega_1}\right).$$
One can verify  that
$g(\b)$ is well-defined  and furthermore depends only on the ideal class
$[\b]$ belonging to in the ideal class group (see
\cite{ramachandra} and \cite{langel}). Let us now describe the main steps
in the proof of the above theorem.

\smallskip
\noindent	
1. In the first step, we need to get an expression for $L(1,f)$.
For this we need  Kronecker's limit
formula. For an ideal  class $\gC$, by Kronecker's limit formula
we have
\begin{equation*}\label{kronecker}
\zeta(s,\gC) =  \sum_{\a \in \gC} {1 \over \N(\a)^s}=
{2\pi \over w\sqrt{|d_K|}}\left(
{ 1\over s-1} + 2\gamma - \log |d_K| - {1 \over 6}\log |g(\gC^{-1})|
\right) + O(s-1), 
\end{equation*}
as $s\to 1^+$. 
Here $d_K$ is the discriminant of $K$ and $w$ is the 
number of roots of unity in $\bO_K$

\smallskip
\noindent
2. In the above expression apart  from $\gamma$, one has the mysterious
number $|g(\gC^{-1})|$. But by CM theory, we know that if
$\gC_1$ and $\gC_2$ are ideal
classes, then $g(\gC_1)/g(\gC_2)$ is an algebraic number lying in
the Hilbert class field of $K$.

\smallskip
\noindent
3. Kronecker's limit formula  gives rise to the following expression for $L(1,f)$
\begin{equation*}
{L(1,f)\over \pi} = {- 1 \over 3w \sqrt{|d_K|}}\sum_{\gC \in \H_K}f(\gC) \log
|g(\gC^{-1})|. 
\end{equation*}
Since $g(\gC^{-1})/g(\gC_0)$ is algebraic for the identity class $\gC_0$
 and  $\rho_f=0$,
we rewrite the above as 
\begin{equation*}
{L(1,f)\over \pi} = {-1 \over 3w \sqrt{|d_K|}}\sum_{\gC \in \H_K}
f(\gC) \log |g(\gC^{-1})/g(\gC_0)|,
\end{equation*}
and hence $L(1,f)/\pi$ is  a linear form in logarithms of algebraic numbers.

\smallskip
\noindent
4. Now we come to the Galois conjugation step. In particular,
one shows that 
 $L(1,f)=0$ implies
that $L(1,f^\sigma)=0$ for any Galois automorphism $\sigma$
of ${\rm Gal}(\bar{\Q}/\Q)$.  This is using Baker's theory
similar to the approach adopted in the  Baker-Birch-Wirsing Theorem.

\smallskip
\noindent
5. In the next step,  one uses the above lemma to 
reduce it to the case when $f$ is actually rational valued.
More precisely, one proves the following: Let $M$ be the  number field
of degree $r$ generated by the values 
of $f$. Then for any basis  $\beta_1, ..., \beta_r$ 
of $M$ over $\Q$ and 
$f(\gC) = \sum_{i=1}^r  f_i(\gC) \beta_i$
with $f_i(\gC) \in \Q$,  we have  $L(1,f)=0$  if and only if
$L(1,f_i)=0$ for $i=1, ..., r$. Thus we may assume without loss of generality  that
our function $f$ is actually rational valued.

\smallskip
\noindent
6. In the penultimate step, it is shown that if $f$ is a rational-valued function and $L(1,f)=0$, then  
$f(\gC)+f(\gC^{-1})=0$ for every ideal class $\gC$.  For this one needs to appeal 
to the Chebotarev density theorem.

\smallskip
\noindent
7. Finally as before, one views these equations as a matrix equation
$DV=0$
where $V$ is the transpose of the row vector $(f(\gC) +f(\gC^{-1}))_{\gC\in H_K}$
and $D$ is the ``Dedekind-Frobenius'' matrix whose $(i,j)$-th
entry is given by 
$\log g(\gC_i^{-1}\gC_j)/g(\gC_i^{-1})$
 with $\gC_i, \gC_j$
running over the elements of the ideal class group.  
The non-vanishing of this determinant is a consequence 
of non-vanishing of  each $L(1,\chi), \chi \neq 1$  and 
consequently $f(\gC)+ f(\gC^{-1})=0$ for all $\gC \neq \gC_0$.
However since the sum
$\sum f(\gC) =0,$ we have $f(\gC_0)=0$ as well,
completing  the proof of  Theorem 13.

\bigskip
\noindent
\section{Linear Independence of L$(k, \chi$)  values with $k > 1$ }
In this final section, let us consider the analogous question for $L(k,\chi)$,
when $k>1$. Here we can also include  the principal character
in our list.

As we saw earlier, the parity of $\chi$ plays a crucial role in 
the context of $L(1,\chi)$, a phenomenon which continues
for $k> 1$. For values in the domain of absolute convergence,
it is worthwhile to  introduce  Hurwitz zeta values as these 
form a natural generating set
for the study of special values of  Dirichlet series 
associated to periodic arithmetic
functions.

For a real number $x$ with $ 0 < x \le 1$ and $s\in \C$ with $\Re(s)>1$,
the Hurwitz zeta function is defined by
\begin{equation*}
\zeta(s,x):=\sum_{n=0}^\infty \frac{1}{(n+x)^s}.
\end{equation*}
This  can be analytically extended  to the entire
complex plane except at $s=1$ where it has a simple pole
with residue one. Note that $\zeta(s,1)=\zeta(s)$.
For any Dirichlet character mod $q$, 
running over arithmetic progressions mod $q$, 
one immediately deduces that 
$$
L(k,\chi) = q^{-k} \sum_{a=1}^{q} \chi(a) \zeta(k,a/q).
$$ 

For $q>1$, let $K_{q}$ be the $\varphi(q)$-th
cyclotomic field. Suppose that $k$ and $\chi$ have the same parity, that is
$\chi(-1) = (-1)^k$. Then the above identity
yields that 
$L(k, \chi)$ in this case is a  $K_{q}$-linear combination of elements of  the following  set
$$
X := \{ \lambda_a ~~:~~1 \le a \le q/2, (a,q) =1\}\,\,\,\,
{\rm where} \,\,\,\,
\lambda_a := \zeta(k, a/q) + (-1)^k \zeta(k, 1-a/q).
$$
Now differentiating the series expansion of  $\pi \cot \pi z$ for $z \notin \Z$, one has
\begin{equation*}
\lambda_a = \zeta(k, a/q) + (-1)^k \zeta(k, 1-a/q) = \frac{(-1)^{k-1}}{(k-1)!} 
~ \frac{d^{k-1}}{dz^{k-1}} (\pi \cot \pi z)|_{(z = a/q)}.
\end{equation*}
On the other hand for $ z \notin \Z$, we have
$$
\frac{d^{k-1}}{dz^{k-1}} (\pi \cot \pi z) 
= \pi^k \sum_{ r, s \ge 0 \atop r + 2s = k} 
\beta_{r,s}^{(k)} \cot^r \pi z 
~(1 + \cot^2 \pi z)^s,
$$
where $\beta_{r,s}^{(k)} \in \Z$.
Since
$i \cot \frac{\pi a}{q} \in \Q(\zeta_q)$,
we see that 
$$
\zeta(k, a/q) + (-1)^k \zeta(k, 1-a/q) = (i \pi)^k \alpha_{a,q},
$$ 
where $\alpha_{a,q} \in \Q(\zeta_q)$.
Thus  when $k$ and $\chi$ have the same parity, we deduce that
$L(k, \chi)$ is an algebraic multiple of $\pi^k$ reminiscent of the fact that $L(1,\chi)$ 
is an algebraic multiple of $\pi$. This also generalises Euler's classical result that
$\zeta(2n)$ is a rational multiple of $\pi^{2n}$.

The following theorem now allows us to settle 
the dimension of the  {\it same-parity} space for prime modulus, namely that
its dimension over $\Q$ is $\frac{p-1}{2}$.
Let $\cot^{(k-1)}(z_0)$ denote the  $(k-1)$-th derivative $\frac{d^{k-1}}{dz^{k-1}} ( \cot z)|_{z = z_0}$.

\begin{thm}\label{lem1}
Let $k>1$ and $q>2$ be positive integers and
 $K$ be a number field such that ${K} \cap \Q(\zeta_q) = \Q$.
Then the set of real numbers
$$
\cot^{(k-1)} (\pi a/q), ~~~~~~~~~~ 1 \le a \le q/2, (a,q) =1
$$
is linearly independent over $K$.
\end{thm}

The above result does seem to be in the spirit of the Baker-Birch-Wirsing theorem
and was proved by Okada \cite{okada}. 
But as noted by Girstmair \cite{Gir}, it is a much simpler result
and we indicate his proof which is short and elegant.
The first point to note is that the numbers
$$
i^k \frac{d^{k-1}}{dz^{k-1}} (\cot  z)|_{(z = \pi a/q)}, ~~~~~~~~~~1 \le a \le q/2,~ ~ (a,q) =1
$$
are Galois conjugates. Now consider any non-trivial $\Q$-linear combination
of the form
$$
\sum_{a} r_a \cot^{(k-1)} (\pi a/q)
= \sum_{a} r_a \left( (-1)^{(k-1)} (k-1)! \frac{q^k}{\pi^k} \sum_{n \in \Z \atop n\equiv a {\rm mod}(q)} n^{-k}\right).
$$
The Galois conjugacy observation allows us to assume that the first coefficient 
$r_1$ has the largest modulus. Thus modulus of the sum in the right hand side of the above identity  is at least
$$
 \frac{(k-1)! q^{k}}{\pi^{k}} |r_1| \left( 1 - \sum_{n=2}^{\infty} n^{-k}\right) > 0
$$
and we are done.  Note that when $k =1$, the matter is
indeed more delicate
as it is linked to the non-vanishing of $L(1,\chi).$

What about the case when $k$ and $\chi$ have opposite parity? Recall for $s=1$,
 the $L(1,\chi)$ values for even characters $\chi$ are linearly independent over $\bar{\Q}$ for all moduli. 
 However for larger integers, the situation is rather bleak with almost no information.
Let us indicate the issue here. When $s=1$, the  Fourier transform approach
 allowed us to express $L(1,\chi)$ as a linear form in logarithms of
 algebraic numbers and thereafter Baker's  theory took over.
An analogous course of action for larger $k$ leads us to 
 linear forms in polylogarithms.

 For an integer $k \geqslant  2$ and complex numbers 
$z \in \C$ with $|z| \leqslant 1$,
the polylogarithm function $Li_k(z)$ is defined as 
$$
Li_k(z) := \displaystyle \sum_{n=1}^{\infty} \frac{z^n}{n^k}.
$$
Then we can deduce the following:
$$
{\rm L}(k,\chi) = \displaystyle \sum_{a=1}^{p} \widehat{\chi}(a) Li_k(\zeta_p^{a}).
$$

However we do not have an analogue of Baker's theorem for  polylogarithms
and hence have  no information when $k$ and $\chi$ have different parity.
We also do not know whether the space generated
by the even and odd $L(k,\chi)$ values intersect trivially.
Presumably, we need deeper results in transcendence
to make any further progress. Finally, there  is a conjecture
of P. Chowla and S. Chowla \cite{CC} on non-vanishing of certain $L(k,f)$
which was later generalised by Milnor \cite{MIL}.
These are deep conjectures which are   linked to the irrationality 
of  numbers of the form $\zeta(2n+1) /\pi^{2n+1}$ as well  as  to  a folklore
conjecture of Zagier on  linear independence of
Multiple zeta values (see
\cite{gmr2},  \cite{tgr} for more details, generalisations  and partial 
results in this direction).

\smallskip
\noindent
{\bf Acknowledgments.} It is our pleasure to thank Sanoli Gun for helpful discussions. The second
 author is a recipient of the MATRICS grant (MTR/2018/000202) from SERB and is 
 thankful to SERB for the support.


\begin{thebibliography}{100}

\bibitem{baker}
A. Baker, {\it Transcendental number theory},
Second edition, Cambridge University Press, Cambridge, 1975.


\bibitem{bbw} 
A. Baker,  B.J. Birch and  E. A Wirsing,
{\it On a problem of Chowla,}
J. Number Theory, {\bf 5} (1973), 224--236. 


\bibitem{BM} 
D. Bertrand and D. Masser, {\it Linear forms in Elliptic Integrals,}
Invent. Math, {\bf 58} (1980), 283--288.


\bibitem{rathbhar} 
A. Bharadwaj and P. Rath, {\it On a question of Alan Baker over number fields,}  Mathematika,
{\bf 66}, Issue 1, 103--111.

\bibitem{tgr} 
T. Chatterjee, S. Gun and P. Rath,
{\it A number field extension of a question of Milnor}, Contemporary Math, {\bf 655},  Amer. Math. Soc. 2015, 15--26.

 
  \bibitem{CC} 
P. Chowla and S. Chowla, {\it On irrational numbers}, 
 Skr. K. Nor. Vidensk. Selsk(Trondheim), {\bf 3}, (1982), 1--5 
(See also S. Chowla, Collected Papers, {\bf 3},
1383--1387, CRM, Montreal, 1999).


 \bibitem{cox}  D. A. Cox, {\it Primes of the form $x^2 + ny^2$}, Willey.
 
  \bibitem{Gir}  K. Girstmair, {\it Letter to the editor}, J. Number Theory, {\bf 23}, 1986, p. 405.
  
  
  \bibitem{GK}
  S. Gun and N. Kandhil,
  {\em On an extension of a question of Baker}, preprint.
  

\bibitem{gmr} S. Gun, M. Ram Murty and P. Rath,
{\it Linear independence of Hurwitz zeta values and a theorem of 
Baker-Birch-Wirsing over number fields,}  Acta Arithmetica, {\bf155}, (2012), no. 3, 297--309.


\bibitem{gmr2}  S. Gun, M. Ram Murty and P. Rath,
{\it On  a conjecture of Chowla and Milnor,}
Canad. J. Math, {\bf 63}, (2011), no. 6, 1328--1344. 



\bibitem{huber} A. Huber and S. M$\ddot{u}$ller-Stach, {\it Periods and Nori Motives}, Springer.

\bibitem{KZ} 
M. Kontsevich and D. Zagier, {\it Periods},  Mathematics 
Unlimited-2001 and Beyond, Springer, (2001),  771--808.

\bibitem{sound} Y. Lamzouri, X. Li and K. Soundararajan,
{\it Conditional bounds for the least quadratic non-residue and related problems,} Math. Comp,
{\bf 84}, (2015), no. 295, 2391--2412.

\bibitem{lang1} S. Lang, {\it Algebraic number theory},
Graduate Texts in Mathematics, 
{\bf 110}, Springer.


\bibitem{Lang1}  S. Lang, {\it Algebraic values of Meromorphic functions I,}
Topology, {\bf 3}, (1965), 183--191.



\bibitem{langel} S. Lang,
{\em Elliptic functions,}
Second edition, Graduate Texts in Mathematics, 
{\bf 112}, Springer.


\bibitem{linnik1}
 Ju. V. Linnik, {\it  On the least prime in an arithmetic progression. I. The basic theorem},
 Rec. Math. [Mat. Sbornik]  N.S., {\bf  15}, (57), (1944), 139--178. 

\bibitem{linnik2}
 Ju. V. Linnik, {\it  
On the least prime in an arithmetic progression. II. The Deuring-Heilbronn phenomenon}, 
Rec. Math. [Mat. Sbornik], N.S., {\bf 15}, (57), (1944), 347--368. 

\bibitem{MIL}
J. Milnor, {\it On polylogarithms, Hurwitz zeta functions,
and their Kubert identities}, Enseignement Math., (2), {\bf 29}, (1983),
no. 3-4, 281--322. 

\bibitem{murty-murty}  M. Ram Murty and V. Kumar Murty,
{\it Transcendental values of class group $L$-functions,}  
Math. Ann., {\bf 351}, (2011), no. 4, 835--855.


\bibitem{ram-kumar-jnt}  M. Ram Murty and V. Kumar Murty,
{\it A problem of Chowla revisited,} 
J. Number Theory, {\bf 131}, (2011), no. 9, 1723--1733.


\bibitem{rath-murty}
M. Ram Murty and P. Rath, {\it Transcendental numbers}, Springer.
 
\bibitem{msaradha}  M. Ram Murty and N. Saradha,
{\it Euler-Lehmer constants and a conjecture of Erd\"os,} 
J. Number Theory, {\bf 130}, (2010), no. 12, 2671--2682.

\bibitem{okada}
T. Okada, {\it On an extension of a theorem of S. Chowla}, Acta Arith., {\bf 38}, (1980/81), 341--345.

\bibitem{ramachandra} 
K. Ramachandra, {\it Some applications of Kronecker's	
limit formulas}, {\sl Ann. of Math} {\bf 80} (1964), no. 2,
104--148.

\bibitem{ram2} 
K. Ramachandra, {\em On the units of cyclotomic fields,}
Acta Arith., {\bf 12}, (1966/67), 165--173.

\bibitem{siegel} C. L. Siegel,
{\em Advanced analytic number theory,}
Second edition,
Tata Institute of Fundamental Research Studies in Mathematics, 
9, 1980.

\bibitem{waldperiods}
M. Waldschmidt, {\it  Transcendence of periods: the state of the art},
 Pure Appl. Math. Q., {\bf 2}, (2006), no. 2, 435--463. 


\bibitem{wash}
L.C. Washington, {\it Introduction to cyclotomic fields,}
 Second edition, Graduate Texts in Mathematics, Springer-Verlag, New York, 1997.
 
 \bibitem{xy} 
 T. Xylouris,
{\it On the least prime in an arithmetic progression and estimates for the zeros of Dirichlet L-functions},
Acta Arith., {\bf 150}, (2011), no. 1, 65--91. 
 
 
\end{thebibliography}
\end{document}